%% file: Residuals-2xn.tex
\documentclass[twoside,12pt, reqno]{amsart}
\usepackage{amsmath,amscd,amsthm,amssymb,amsxtra,latexsym,epsfig,epic,graphics}
\usepackage[matrix,arrow,curve]{xy}
\usepackage{graphicx}
\usepackage{diagrams}
\usepackage{tikz,color}  
\usepackage[shortlabels]{enumitem}
\usepackage{showlabels}
\input preamble.tex

\showlabels{}

\def\FF{{\mathbb F}}

\def\fix#1{{\bf ***Fix:} #1 {\bf ***}}

\DeclareMathOperator{\rH}{{\rm H}}

\def\rH{{\rm H}}

\def\RR{{\mathcal R}}

\def\mm{{\frak m}}

\def\RR{{\mathcal R}}

\def\grade{{\rm grade}}

\newarrow{Iso} -----

\def\m{{\frak m}}

\def\lbracket{{[\kern-1.5pt[}}
\def\rbracket{{]\kern-1.5pt]}}

\def\ite#1{{\ref{#1}}}
\newtheorem{proposition-definition}[theorem]{Proposition-Definition}

\makeatletter
\def\Ddots{\mathinner{\mkern1mu\raise\p@
\vbox{\kern7\p@\hbox{.}}\mkern2mu
\raise4\p@\hbox{.}\mkern2mu\raise7\p@\hbox{.}\mkern1mu}}
\makeatother

\usepackage{times}
\newdimen\x \x=12pt

\usepackage{color}

\usepackage[breaklinks,bookmarksopen,bookmarksnumbered,urlcolor=blue]{hyperref}
\hypersetup{colorlinks=true,backref=true,citecolor=blue}

\date{\today}
\author{David Eisenbud and Bernd Ulrich
}

\title{Residual Intersections of $2\times n$ Determinantal Ideals}
\begin{document}

\begin{abstract} 
Schemes defined by residual intersections have been extensively studied in the case when they are Cohen-Macaulay, but this is a very restrictive condition. In this paper we make the first study of a class of natural examples far from satisfying this condition, the rank 1 loci of generic $2\times n$ matrices. Here we compute their depths and many other properties. These computations require a number of novel tools.
\end{abstract}

\maketitle
\section*{Introduction}
 
\let\thefootnote\relax\footnote{
\noindent
AMS Subject Classification:
Primary 13C40, 13H10, 14M06, 14M10, 14M12;
Secondary 13D02, 13N05, 13P20.

\noindent
We are grateful to the NSF for partial support through awards 2001649 and 1802383, and to MSRI for hosting us during some of the work on this paper.}

If $I\subset S$ is an ideal in a Noetherian ring, and $J\subset I$ is an ideal generated by $s$ elements, then the ideal $K := J:I$ is called an \emph{$s$-residual intersection} of $I$
if $\codim K \geq s$, and a \emph{geometric $s$-residual intersection} of $I$ if, in addition, the codimension of $I+K$ is at least $s+1$. If $S$ is Gorenstein, $I$ is unmixed, and $s=\codim I$, and $K\neq S$ then $K$ is said to be a \emph{direct link} of $I$, and in this case the properties of $K$ are very tightly
related to those of $I$; for example, if $S/I$ is Cohen-Macaulay, then so is $S/K$, and if the link is geometric then the canonical module of $S/K$ is $(I+K)/K$. 
However, when $s>\codim I$, a case of great geometric interest, the connection is not nearly
as simple. For example, much stronger conditions than $S/I$ being Cohen-Macaulay are necessary for $S/K$ to be Cohen-Macaulay.

Residual intersections have been studied from an enumerative point of view since the work of Chasles resolving the Steiner problem ("How many plane conics are tangent to 5 general plane conics") in 1864, and reappeared independently in the work of Brill, Noether, Macaulay and others on the Riemann-Roch and Cayley-Bacharach theorems in the 1880s and 1890s, as well as in Halphen's famous 1882 classification of projective space-curves of low degree. The enumerative theory was put into rigorous modern form by Fulton and MacPherson, and codified in Fulton's monumental 1984 work \cite{F}.
 
Scheme-theoretic (equivalently, ideal-theoretic) results on links have been known (in special cases) since Macaulay's work in 1913 \cite{M}, 
and were independently and more generally codified in 1974 by Peskine and Szpiro in~\cite{PS}. The modern study of the ideals defining higher residual intersections was initiated by
Artin and Nagata \cite{AN}, focusing on the question of when the residual intersection is Cohen-Macaulay, and this question has been studied extensively by Huneke, Ulrich and others in \cite{H2}, \cite{HVV}, \cite{HU1}, \cite{U}, \cite{CEU}, \cite{EU1},  \cite{H}, \cite{HN}, \cite{CNT},\cite{EHU}. A special feature of this development is that under good circumstances not only is $R = S/K$ Cohen-Macaulay, but
the canonical module is $(IR)^{s-\codim I+1}$, and $IR$ generates the class group of $R$.

These properties hold, for example, when $I$ is the ideal of minors of a sufficiently general $(n-1)\times n$ matrix; but they fail already
for ideals as simple as the ideal of $2\times 2$ minors of a $2\times 4$ matrix, and for many other simple ideals.

\begin{figure}\label{Fig1}
\caption{Depth of $(IR_i)^j$}
\centerline{
\includegraphics[width=200mm]{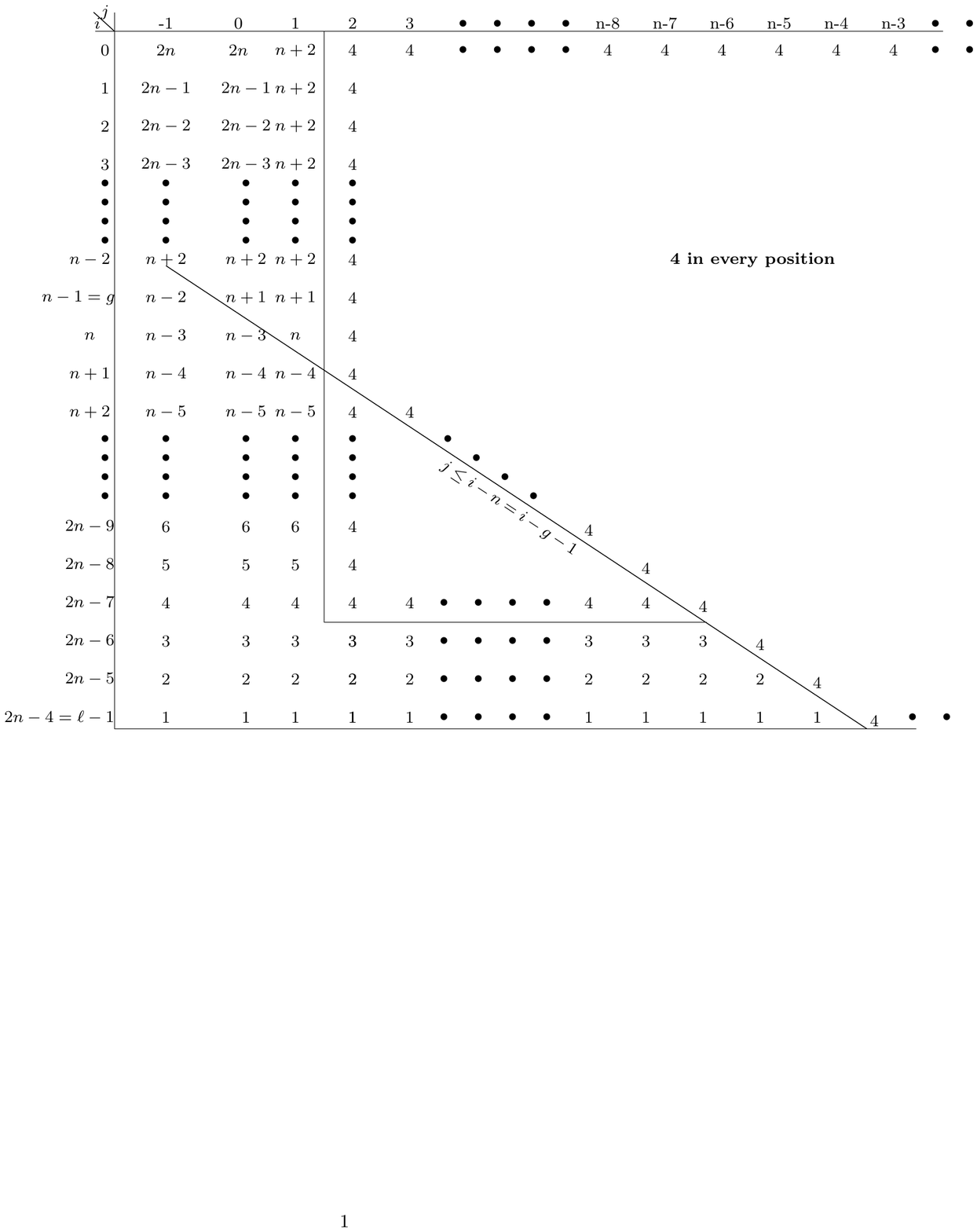} \qquad\qquad\qquad \qquad\qquad\qquad \qquad
}
\end{figure}

In this paper we initiate a detailed study of the ``first" case when the strong conditions fail, namely the ideal $I$ of $2\times 2$ minors of a generic $2\times n$ matrix,
 $$
\begin{pmatrix}
 x_{1,1}&x_{1,2}&\dots&x_{1,n}\\
x_{2,1}&x_{2,2}&\dots&x_{2,n}
\end{pmatrix}
 $$
for arbitrary $n\geq 4$ (the case $n=3$ being easy and special).  We compute the depths of the sufficiently general $i$-residual intersections $R_i$ for $i\leq \ell(I)-1$ (where $\ell(I)$ is the \emph{analytic spread}  of $I$, the dimension of the Grassmannian $\bold {G}(2,n)$ plus 1), the case $i\geq \ell(I)$ being trivial,
and the depths of  the ideals
$(IR_i)^j$ for $-1\leq j$ (Theorems~\ref{new} and~\ref{main}). A table of the depths
is given in Figure~\ref{Fig1}.

The importance of the ideals $(IR_i)^j$ is suggested by Theorem~\ref{class group} which says that under additional genericity assumptions, $IR_i$ generates
the class group of $R_i$. Their properties are also needed in our computation of the depths of the residual intersections themselves.

Since the proofs of our main results involve a complicated induction, it seems worth pointing out some of the main new ideas of the paper.
To understand the situation, recall that two technical conditions on $I$ have played a central role in the theory of residual intersections, culminating in the paper~\cite{U}.
 The first, $G_s$, ensures that  $s$-residual intersection $K = (a_1,\dots, a_s):I$ will have codimension $s$ whenever the choice of $a_1,\dots, a_s$ is sufficiently general.
 This condition is trivially satisfied for the rank 1 loci of $2\times n$ matrices.
 
 The second condition is much more serious, and fails dramatically for our class of ideals. It is a condition on the depths of the ``thickenings'' $S/I^j$, and may be stated as
$$
\Ext_S^k(S/I^j, S) =0 \quad \hbox{ for all $j\leq s-g+1$ and $k \geq g+j \, ,$ }
$$
where $g = \codim I$. 

The first advance that enables our computations is~\cite[Corollary 4.2]{CEU}, where we showed that the $s$-residual intersection of 
the ideal $I$ satisfies Serre's condition $S_2$ if the much weaker vanishing condition
$$
\Ext^{g+j}_{S}(S/I^{j}, S) =0 \quad \hbox{ for $1\leq j\leq s-g +1$}
$$ 
holds.  By~\cite[Theorem 4.3]{RWW} this condition is satisfied in our case.

We use this to compute the depths of the ideals $(IR_i)^j$ for large $j$, but need a form of duality to compute even the depth of $R_i$ itself. For this we 
need to prove that the canonical module, even in this non-Cohen-Macaulay case, is given by the formula $\omega_R= (IR)^{s-g+1}$. The usual proofs do no work in our case, so instead we use our \cite[Theorem 4.1]{EU1}, which gives a map $(IR)^{s-g+1}\to \omega_R$ defined under very weak hypotheses. In our case these are satisfied,
and we are able to show that the map is an isomorphism under much weaker conditions than were known before.

Other important inputs to our proof include the representation theoretic analysis of free resolutions of powers of $I$ in \cite{RWW} and a surprising equality
between certain betti numbers of $S/I^j$ and $S/I^{j+2}$ whose proof was shown to us by Raicu after we had discovered the phenomenon experimentally using
\cite{M2}.

It may be asked why we treat the case of a $2\times n$ matrix but not $3\times n$ or larger matrices. The fact is that residual intersections in most of these cases
simply do not exist: the  Artin-Nagata condition $G_s$, which is an assumption on the local number of generators at various primes,
 is rarely satisfied in determinantal cases beyond $2\times n$.

 In \cite{EHU} we made
an analysis of the $(\ell(I)-1)$-residual intersections of $I_2(X)$, and various generalizations, showing
in particular that for large $j$ the module $(IR_{2n-4})^j$ is Cohen-Macaulay (a special case
of Theorem~\ref{main}\ref{depth} below). In Section~\ref{asymp sect} we give a general bound on the asymptotic depth of the  powers of an arbitrary ideal $I$ modulo its $s$-residual intersection, showing that except for $s=\ell(I)-1$, the high powers are never  Cohen-Macaulay modules.

In~\cite{EU2}, in progress, we study the Rees algebras of the ideals $IR_i\subset R_i$.
\subsection*{Note on computation} The research that went into this paper was heavily dependent on computations using Macaulay2, without which we would not have guessed the
rather unintuitive form of the result of the main Theorem. In making these computations, we found it useful to use a rather sparsely generated
sequence of ideals $J$---this allowed us to make the difficult residual ideal computations in many more cases. The idea is that
there is a minimal reduction of the ideal of $m\times m$ minors of a generic $m\times n$ matrix
 (or its initial ideal) consisting of the $\ell = m(n-m)+1$ sums of these minors whose column numbers
add to a given number (see \cite[Theorem 6.3]{DEP}). For example, for a $2\times 5$ matrix with columns $1,\dots,5$ the sums of minors may be represented by:
$$
(1,2), (1,3), (1,4)+(2,3), (1,5)+(2,4), (2,5)+(3,4), (3,5), (4,5).
$$
The initial subsequence of $i$ elements of this sequence, in the case $m=2$, generates an ideal $J_i$
that makes the computation of the residual intersection $J_i:I$ relatively quick.
\medskip

\section{Depths of Residual Intersections}
We use the following notation throughout this paper:
Let $$
 X = \begin {pmatrix} x_{1,1}&\dots& x_{1,n}\\
  x_{2,1}&\dots&x_{2,n}
  \end {pmatrix}
  $$
  be a generic $2\times n$ matrix over the ring 
 $S = k[x_{1,1},\dots,x_{2,n}]$, where $k$ is a field of characteristic 0 and $n\geq 4$ (the case $n=3$ being easy and special). Write $\m$ for the maximal homogeneous 
 ideal of $S$, and let $I= I_2(X)$ be the ideal of $2\times 2$ minors of $X$. The codimension and analytic spread of $I$ are
 $g := n-1, \ell := 2n-3$ respectively. Let $a_1,\dots,a_{\ell}$ be quadrics in $I$ that generate a minimal reduction. Set
 $J_i = (a_1,\dots, a_i)\subset I$ and $R_i = S/(J_i:I)$
 and assume that $K_i:=J_i:I$ is a geometric $i$-residual intersection for $i \leq \ell-1$, which will be the case
 for a general choice of the $a_i$. 
We henceforth assume that $0\leq i\leq \ell-1$ because the $\ell$-residual intersection would be primary to the homogeneous maximal ideal. We set $(IR_i)^0:=R_i$,
and if $j<0$ we write $(IR_i)^{j}$ for the inverse ideal $((IR_i)^{-j})^{-1} = R_i:_{{\rm Quot}(R_i)}(IR_i)^{-j}$ . 

\smallskip

\begin{theorem}\label{new} With hypotheses as above:

 \begin{enumerate}[$($a$)$]
\item \label{CM,Bb} The rings $R_i$ are Cohen-Macaulay on the punctured spectrum. They are Cohen-Macaulay if $i\leq g$ and, at least for $n\geq 5$,  Buchsbaum if $i=g+1$. If  $i\geq g+1$ then  
$$
\depth R_i=2n-3-i = \dim R_i-3,
$$ 
and the local cohomology of $R_i$ is nonzero only in cohomological degrees 
$
   \depth R_i   \hbox{ and } \dim R_i .
$
 
\item\label{unmixed} If $g+1\leq i\leq \ell -2$ the ideal $IR_i$ is unmixed of codimension 1, $\depth R_i/(IR_i)\geq \depth R_i$,
and $\depth (IR_i)^{-1}/R_i\geq \depth R_i$.

\item \label{can}For $i\geq 0$ the canonical module of $R_i$ is 
$$
 \omega_{R_i} = (IR_i)^{i-g+1}(2(i-g-1)),
 $$ 
noting that $(IR_i)^{i-g+1}=R_i$ if $i \leq g-1$.

\item\label{regularity}
For $i\leq g-1$ the ring $R_i$ has regularity $i,$ for $i = g$ it has regularity $i-2,$ and for $i = g+1 \geq 5$ it has regularity
 $i-3.$ 
\end{enumerate}
\end{theorem}

\smallskip

\begin{conjecture}
For $i\geq g+1 \geq 5$ the regularity of $R_i$ is $i-3.$
\end{conjecture}

\smallskip

\begin{theorem}\label{main}
With hypotheses as above and $j\geq -1:$
\begin{enumerate}[$($a$)$]\label{iso of M-main}

\item\label{iso} $(IR_i)^j \cong  I^j/J_iI^{j-1}$ for $j>0.$

\item\label{pspec} The modules $(IR_i)^j$ are locally Cohen-Macaulay on the punctured spectrum for
$ j\leq i-g+2.$

\item\label{loCo} The local cohomology of the modules  $IR_i$, $R_i$ and $(IR_i)^{-1}$ are 0 except at the depth, given below, and dimension, $2n-i.$

\item\label{loco duality} If $j\leq i-g+2$ and one of the conditions
\begin{enumerate}[(i)]
\item  $i\leq \ell-2;$ or 
\item $i=\ell-1$ and $j\geq i-g,$
\end{enumerate}
is satisfied, then the multiplication maps $(IR_i)^j\otimes (IR_i)^{(i-g+1)-j} \to (IR_i)^{i-g+1}$ induce duality isomorphisms
$$
\Hom((IR_i)^{(i-g+1)-j} , \omega_{R_i}) \cong (IR_i)^j(2(i-g-1)),\qquad (1)
$$
and thus $(IR_i)^j$ satisfies
Serre's condition $S_2$.
Furthermore, if $2\leq p\leq \dim R_i-1$, then
$$
\rH^{\dim R_i+1-p}_\mm((IR_i)^{(i-g+1)-j} )^\vee \cong \rH^{p}_\mm ((IR_i)^j)(2(i-g-1)),\quad (2)
$$
where $-^{\vee}$ denotes $k$-dual.

\item\label{depth}
$$
\depth (IR_i)^j = 
\begin{cases}
 \dim R_i& \hbox{if $j\leq 0$ and $i\leq g$ }\hfill (3)\\
\dim R_i -3 & \hbox{if $-1\leq j\leq \min\{1, i-g-1\}$} \hfill\quad (4)\\
 n+2 & \hbox{if $j=1$ and $i\leq g-1$}\hfill (5)\\
\dim R_i &\hbox{if $j=1$ and $g\leq i\leq g+1$ }\hfill (6)\\
 \min\{\dim R_i -3,4\} & \hbox{if $2\leq j\leq i-g-1$}\hfill (7)\\
 4 &\hbox{otherwise.}\hfill (8)
\end{cases}
$$
\end{enumerate}
\end{theorem}

\bigskip
\medskip

\noindent {\bf Remark:} In the cases $g=4, i = 4$ (linkage) and the case
$g=3, i = 2$ the ideal $(IR_i)^3$ can have an associated
prime of codimension 2 in $R_i$ (by computation in an example.) Thus
the $S_2$ conclusion of Theorem~\ref{main}(d) can fail outside the given range.

As a consequence, we see that $M_{i,j}$ is a maximal Cohen-Macaulay $R_i$-module if and only if:
\begin{itemize}
 \item $j=0$ and $i\leq g$; or
 \item $j=1$ and $g-1\leq i\leq g+1$; or
 \item $j\geq n-3$ and $i = \ell-1$.
\end{itemize}

In the proofs we will use the following statements, which are essentially in the literature:

Recall that an ideal $I$ is called weakly $s$-residually $S_2$ if,
for any geometric $i$-residual intersection $K$ of $I$ with $i\leq s$,
the ring $S/K$ satisfies Serre's condition $S_2$.

Also, recall that an ideal $J$ is said to be of linear type if the natural surjections
$$
\Sym_j(J )\rOnto  J^j 
$$
are isomorphisms.
 \vfill\break
 \goodbreak
 
  \begin{proposition}\label{known facts}
Let $S$ be any Cohen-Macaulay ring,  and $J := (a_1, \dots, a_s) \subsetneq I \subset S$
any ideals. Set $J_i := (a_1,\dots, a_i)$. 
Suppose that $I$ is weakly $(s-2)$-residually $S_2$ and that for $i<s$ the ideal
$K_i = I:J_i$ is a geometric $i$-residual intersection of $I$. 
For $i< s$ we have:

\begin{enumerate}[label=(\alph*), ref = (\alph*)]
\item \label{colon condition} $J_i:a_{i+1} = J_i:I$.
\item\label{intersection} $(J_i:I)\cap I = J_i$.
 \item \label{Ki unmixed} $K_i$ is unmixed of codimension $i$.
 \item \label{nzd} The element $a_{i+1}$ is a nonzerodivisor on $S/K_i$.
  \item \label{intersection2} $J_i\cap J^j  = J_i J^{j-1}$ if $j>0$.
 \item \label{linear type} $J/J_i \subset S/J_i$ is an ideal of linear type.
 \end{enumerate}
\end{proposition}
\begin{proof}

From the assumptions we see that the ideal $I$ satisfies the condition $G_s$. 
 Because $I$ is weakly $(s-2)$-residually $S_{2}$, it is $(s-1)$-parsimonious and $(s-1)$ thrifty in the sense of \cite[Proposition 3.1]{CEU}. This means that for every $0\leq k\leq s-1$ we have
\begin{align*}
 J_k:a_{k+1} &= J_k :I  &\hbox{(parsimonious)}\\
 (J_k:I)\cap I &= J_k &\hbox{(thrifty)}.
\end{align*}
In particular, this proves parts~\ref{colon condition} and \ref{intersection}.

 Note that $K_i$ is a proper ideal because $J\neq I$. Thus we may 
 apply~\cite[Proposition 3.3(a)]{CEU} and conclude that $K_i$ is unmixed of codimension exactly $i$. Also, \cite[Proposition 3.3(b)]{CEU}
 is the statement \ref{nzd}.  
 
We postpone the proof of part~\ref{intersection2} to prove part~\ref{linear type}. Combining parsimony and thrift and using $J\subset I$ we get 
$$
(J_k:a_{k+1})\cap J = J_k.
$$
For $0\leq i\leq s-1$ we set $S_i = S/J_i$. For $0\leq i\leq k\leq s-1$ we have
$$
(J_kS_i:a_{k+1}S_i ) \cap JS_i = J_k S_i;
$$
because $J_i\subset J_k$ and $J_i\subset J$.
This means that the images of $a_{i+1},\dots,a_s$ in $S_i$ form a relative $*$-regular sequence as defined in \cite{HSV}. It
follows by \cite[Theorem 5.6]{HSV} that the $\mathcal M$-complex of $JS_i$ is acyclic. Therefore 
by \cite[Corollary 2.2 and Theorem 2.3]{HSV} $JS_i$ is of linear type. (See also \cite[Theorem 3.1]{H1} or \cite[Theorem 3.15]{V}.)

For part \ref{intersection2}, we assume $j>0$ and consider the diagram
$$
\begin{diagram}
\frac{\Sym_j(J) }{\Sym_1(J_i)\Sym_{j-1}J}& \rOnto& \frac{J^j}{J_i J^{j-1}}\\
\dTo^{\mu} &&\dOnto_\beta\\
\Sym_j(\frac{J}{J_i}) &\rOnto^\alpha & \frac{J^j}{J_i\cap J^j}.
\end{diagram}
$$
Since $\mu$ is an isomorphism and $\alpha$ is an isomorphism by part \ref{linear type}, the map $\beta$ is injective, and thus an isomorphism, as required.
\end{proof}

We now return to the notation and assumptions in the first paragraph of this section.

\begin{corollary}\label{residually S2}
 In characteristic 0, $I$ is weakly $(\ell-2)$-residually
$S_2$, and thus the conclusions of Proposition~\ref{known facts} hold for the ideal $I$ with
$s =  2n-3$. In particular, the dimension of $R_i$ is $2n-i$, as asserted
in Theorem~\ref{new}\ref{CM,Bb}.
\end{corollary}

\begin{proof}Recall that $\ell := \ell(I) = 2n-3.$
By~\cite[Theorem 4.3]{RWW}, $\Ext^{n+j-1}_{S}(S/I^{j}, S) =0$ for $2\leq j\leq n-3 = (\ell-2)-\codim I +1$ (this is where we require characteristic 0). The same formula holds for $j=1$ because $S/I$ is Cohen-Macaulay of codimension $n-1$. By~\cite[Corollary 4.2]{CEU}, this implies that $I$ is $(\ell-2)$-residually $S_{2}$,
verifying the hypotheses of Proposition~\ref{known facts}.
\end{proof}

\begin{proof}[Proof of Theorems~\ref{new} and \ref{main}]

We begin by proving Theorem~\ref{main}\ref{pspec}. The first statement of
Theorem~\ref{new}(a) is the case $j=0$ of that statement.

 The ideal $I$ is a complete intersection
on the punctured spectrum. The case $j=0$ is the statement that
$R_i$ is Cohen-Macaulay on the punctured spectrum, proved in \cite[Theorem 3.1]{H}. 
The case $j>0$ follows from \cite[Corollary 3.10]{HU2}.
To handle $j=-1$ we can, on the punctured spectrum, write 
$$
(IR_i)^{-1}\cong \Hom_{R_i}(IR_i, R_i) \cong \Hom_{R_i}(IR_i \otimes \omega_{R_i}, \omega_{R_i})
 \cong 
 \Hom_{R_i}((IR_i)^{i-g+2}, \omega_{R_i})
 $$
 where the first isomorphism follows from Proposition~\ref{known facts}\ref{nzd}, the second isomorphism
 uses the fact that $R_i$ is Cohen-Macaulay on the punctured spectrum, and
 the third isomorphism uses the computation of the canonical module \cite[Proposition 2.3]{HU1} and the fact that
 the kernel of the surjection $IR_i \otimes \omega_{R_i} \to (IR_i)^{i-g+2}$ is torsion,  because
 $R_i$ is generically a complete intersection, completing the proof of Theorem~\ref{main}\ref{pspec}.
 
Set
\begin{align*}
 M_{i,j} &:= (I^j/J_iI^{j-1})(2j) \hbox{ if } j>0,\\ M_{i,j} &:= (IR_i)^j(2j) \hbox{ if } j\leq 0.
\end{align*}
 We will prove the statements in both Theorems together with all of the following items, by induction on $i$.

\begin{enumerate}[label=(\Alph*)]
 \item  \label{local coh} $\rH^p_\mm (M_{i,j}) = 0$ if $j\geq i-g$ and $5\leq p \leq \min\{n-j+2,\ \dim R_i -1\}$.

\item \label{depth for M} The depth of each $M_{i,j}$ is at least what is given as the depth
of $(IR_i)^j$ in  Theorem~\ref{new}\ref{CM,Bb} and Theorem~\ref{main}\ref{depth}.
 
 \item\label{iso of M} The natural maps $M_{i,j} \to (I R_i)^j(2j)$ are isomorphisms.
 
 \item \label{hooks} The maps 
 $$M_{i,j} \rTo^{a_{i+1}} M_{i,j+1}$$
 are monomorphisms, so that the natural sequences
 $$
 \leqno{(*_{i,j})} \qquad\qquad 0\rTo M_{i,j} \rTo^{a_{i+1}} M_{i,j+1} \rTo M_{i+1,j+1} \rTo 0
 $$
 are exact.
 
 \item  \label{linearity}Let $F^{i,j}$ be the minimal graded free resolution of $M_{i,j}$ over $S$. 
 
\begin{enumerate}
 
 \item If $j\geq i$ then $F^{i,j}$ is linear.
 
 \item If $j\geq i-g$  then $F^{i,j}_{k}$ is generated in degree $k$ for $k\geq i+1$.
 
 \item If $j\leq  i-g-1$ then $F^{i,j}_{k}$ is generated in degree $k$ for $k\geq i+4$.
\end{enumerate}

\end{enumerate}

\def\it#1{{(\ref{#1})}}

\smallskip

We first deal with the case $i = 0$, noting that $M_{0,j} = I^j(2j)$. Theorem~\ref{new} and parts (a) and (d) of Theorem~\ref{main} are obvious for $i=0$.
Theorem~\ref{main}(c), and item \ref{local coh} in the cases $-1\leq j\leq 1$, follow because $S$ and $S/I$ are Cohen-Macaulay. The 
cases $j\geq 2$ of item \ref{local coh} are a consequence of \cite[Theorem 4.3]{RWW}. 
Items \ref{depth for M} and \ref{linearity} and Theorem~\ref{main}(e) follow from \cite[Theorem 5.4 and the beginning of its
proof]{ABW}. Items \ref{iso of M} and \ref{hooks} are obvious for $i=0$. 
\medbreak
Now assume $i>0$ and that the assertions all hold for $i-1$. 

\bigbreak
To prove item \ref{local coh} in case $j=-1$ note that $i$ must be $\leq g-1$, so
$R_i$ is  a complete intersection of codimension $i$. If $i\leq g-2$ then $IR_i$ has codimension $\geq 2$, so $M_{i,-1} = (IR_i)^{-1}(-2) = R_i(-2)$. On the other hand if $i = g-1$, then $IR_i$ is a height 1 ideal and $R_i/(IR_i) = S/I$, is Cohen-Macaulay, so
$IR_i$ is a maximal Cohen-Macaulay module over the Gorenstein ring $R_i$.
Thus its dual, which is $M_{i,-1}$ up to a twist, is again a maximal Cohen-Macaulay module. The local cohomology vanishing of part 
\ref{local coh} follows. If $j\geq 0$ the statement of \ite{local coh} follows from the long exact sequence in cohomology applied to the exact sequence $(*_{i-1,j-1})$ of \ite{hooks} and item \ite{local coh} for $i-1$. This also proves Theorem~\ref{main}\ref{loCo} for $(IR_i)^j$ in case $i \leq g+j$.

\bigbreak

We next consider item \ite{depth for M} for $i$, postponing the case $j=-1$.

First suppose $0\leq j\leq 1$. By the induction on $i$,  $\depth M_{i-1,j-1}$ is greater than
the claimed $\depth M_{i,j}$ and $\depth M_{i-1,j}$ is at least as large as the claimed depth, so the result follows from \ite{hooks}. 

Next, suppose that $2\leq j\leq i-g-1$. The projective dimension of $M_{i-1,j}$ is $\leq 2n-\min \{2n-3-(i-1),4 \}=\max\{i+2,2n-4\}$ and by \ite{linearity} the minimal resolution of $M_{i-1,j-1}$ is
linear from step $(i-1)+4 = i+3$ on. 
Using Lemma~\ref{linearitylemma}(a) we see that the projective dimension
of $M_{i,j}$ is at most $\max\{2n-4, i+3\} = 2n - \min\{\ell -i, 4\}$, as required.

If $j\geq \max\{2,i-g\}$ then, by induction $\depth M_{i-1, j} \geq 4$, so the projective dimension of $M_{i-1,j}$ is $\leq 2n-4$. By \ite{linearity} the minimal resolution of
$M_{i-1,j-1}$ is linear from step $(i-1)+1 = i\leq 2n-4$ on. Now again Lemma~\ref{linearitylemma}(a) implies that $\pd M_{i,j}\leq 2n-4$, as required.

\bigbreak
We now  prove  items \ite{iso of M}, \ite{hooks} for $i$ and for all $j$.

For $j\leq 0$, item \ite{iso of M} is trivial. The case $j=-1$ of \ref{hooks}
follows because
$a_{i+1}$ is a nonzerodivisor on $R_i$ by Proposition~\ref{known facts}\ref{nzd}. The case $j=0$ follows by Proposition~\ref{known facts}\ref{colon condition}, which implies that
$$
S/(J_{i}:I) \rTo^{a_{i+1}} I/J_{i}
$$
is a monomorphism. 

Now assume that $j>0$.  We claim that 
$$
 \leqno{(**)} \qquad\qquad (J_{i}: a_{i+1})\cap I^{j}  = (J_i:I)\cap I^j = J_iI^{j-1}.
$$
 The first equality follows from Proposition~\ref{known facts}\ref{colon condition}. Since $j>0$,
 $$
(J_i:I)\cap I^j = (J_i:I)\cap I\cap I^j = J_i\cap I^j,
$$
where the last equality follows from Proposition~\ref{known facts}\ref{intersection}.

Thus to complete the proof of $(**)$ we must show that 
$$
 \leqno{(***)} \qquad\qquad J_{i}\cap I^{j} = J_iI^{j-1}.
$$
We first prove equality locally on the punctured spectrum. Let $P\neq \mm$ be a homogeneous 
nonmaximal prime ideal containing $I$. Note that $I_{P}$ is a complete intersection. Because $J:=J_\ell$ is a reduction
of $I$ we see that $I_{P} = J_{P}$. Now Proposition~\ref{known facts}\ref{intersection2} implies that
$(J_{i}\cap I^{j})_{P} = (J_{i}I^{j-1})_{P}$.

By definition $J_i\cap I^j/J_iI^{j-1} \subset M_{i,j}$ and $M_{i,j}$ has positive depth by \ite{depth for M} with $j\geq0$. Since
$J_i\cap I^j/J_iI^{j-1}$ is zero on the punctured spectrum, it must have finite length and thus is 0.

The second equality of $(**)$ implies
\ref{iso of M} for $j>0$. The equality of the first and third terms in $(**)$ implies 
$$
(J_{i}I^{j}:a_{i+1})\cap I^{j}= J_{i}I^{j-1},
$$
which gives \ref{hooks} for $j>0$. This concludes the proof of \ite{iso of M} and \ite{hooks} for $i$, and with it Theorem~\ref{main}\ref{iso}.

\bigbreak
We next prove Theorem~\ref{new}\ref{can} and
Theorem~\ref{main}\ref{loco duality}. If $i\leq g-1$ the ring $R_i$ is a complete intersection, and the results are easy,
so we may suppose $i\geq g$. 
By \cite[Theorem 4.1]{EU1} there is a map $\mu: M_{i,i-g+1}(-4) \to \omega_{R_i}$, which is an isomorphism on the punctured spectrum because $I$ is a complete
intersection there. Notice that $i-g+1\geq \max\{1, i-g\}$, and by \ite{depth for M} with $j>0$, the depth of $M_{i,i-g+1}$ is $\geq 2$. Thus $\mu$ is
an isomorphism. 

Furthermore, the multiplication maps $(IR_i)^j\otimes (IR_i)^{(i-g+1)-j} \to (IR_i)^{i-g+1}$ induce maps
$$
(IR_i)^{j}(-4) \to \Hom((IR_i)^{(i-g+1)-j}, \omega_{R_i}).
$$ 
These maps become isomorphisms locally at each prime $P$ on the punctured spectrum:
Since $\omega_{R_i} = (IR_i)^{i-g+1}$ up to shift, the map is the natural isomorphism at
$P$ as long as $(IR_i)_P$ is principal and generated by a nonzerodivisor. 
By Theorem~\ref{main}(b)
it suffices, for the proof on the punctured spectrum, to prove the isomorphism  for primes $P$ containing $IR_i$ and of codimension 1 
in $R_i$. But $I+K_{i+1} = I+(J_{i+1}:I)$ has codimension $\geq i+2.$ Thus locally at $P$, $IR_i$ is generated by $a_{i+1}$, which is a nonzerodivisor by
Proposition~\ref{known facts}\ref{nzd}. This completes the proof of the fact that the maps above are isomorphisms locally 
on the punctured spectrum.

Hence, to prove that these maps are isomorphisms globally it suffices to show that $\depth (IR_i)^j \geq 2$. 
Suppose that $i\leq \ell-2$. If $j\geq 0$ then the assertion follows from \ite{depth for M} for $j \geq 0$.
As, in particular, $\depth R_i \geq 2$ by \ite{depth for M} and $IR_i$ contains a nonzerodivisor on $R_i$ by
Proposition~\ref{known facts}\ref{nzd}, we conclude that the inverse ideal $(IR_i)^{-1}$ has depth at least 2, proving
the assertion for $j=-1.$ Finally, if $i = \ell-1$ and $j\geq i-g$, we may again apply \ite{depth for M} since $j \geq 0.$

Theorem~\ref{main}\ref{loco duality}(2) follows from Lemma~\ref{duality lemma} because
$R_i$ is equidimensional of dimension $d-i$ by Proposition~\ref{known facts}\ref{Ki unmixed},  $(IR_i)^j$ contains a nonzerodivisor on $R_i$ by Proposition~\ref{known facts}\ref{nzd}, 
and this module is Cohen-Macaulay locally on the punctured spectrum by Theorem~\ref{main}\ref{pspec}.
\bigbreak

We can now prove \ite{depth for M} for $j=-1$ and Theorem~\ref{main}\ref{loCo}. If $i\leq g-1$ and $j=-1$
then $M_{i,j}$ is a maximal Cohen-Macaulay $R_i$-module by the argument at the beginning of the proof of \ite{local coh} above.
If $i\leq g-1$ and $0\leq j \leq1$, Theorem~\ref{main}\ref{loCo} follows from the exact sequence
$$
0\to I/J_i \to S/J_i \to S/I\to 0 \, .
$$

Now suppose that $i\geq g$. As $i-g+1-j \geq i-g \geq 0$, we may use \ite{local coh} and \ite{depth for M} in the case of
exponent $i-g+1-j$ to say that $\rH^p_{\mm}((IR_i)^{i-g+1-j})=0$ for $5 
\leq p \leq 2n-i-1$ and that $\depth (IR_i)^{i-g+1-j} \geq 4$. Thus $(IR_i)^{i-g+1-j}$
has at most  one nonvanishing intermediate local cohomology module, $\rH_\mm^4((IR_i)^{i-g+1-j})$.  Therefore
by the duality statement, Theorem~\ref{main}\ref{loco duality}, equation (2), 
 the only possible nonvanishing intermediate local cohomology modules of $(IR_i)^j$ are $\rH_\mm^q((IR_i)^{j})$ for
 $q$ equal to $0$, $1$, $\dim R_i -3$. The first two possibilities, $q=0$ and $q=1$, are ruled out by 
 Theorem~\ref{main}\ref{loco duality} when $i \leq \ell -2=2n-5$. If $i=\ell-1=2n-4$, then $1= \dim R_i-3$ and $\depth (IR_i)^j \geq 1$. The last inequality
 follows from \ite{depth for M} if $j\geq 0$ and holds for $j=-1$ because $R_i$ has depth $\geq 1$.
 This concludes the proof  of \ite{depth for M} for $j=-1$ and  of Theorem~\ref{main}\ref{loCo}, and this also proves
 the last statement of
 Theorem~\ref{new}\ref{CM,Bb}, for $i$.

\bigbreak

Continuing the induction, we next prove \ite{linearity}. First suppose that $j=-1$. If $i\leq g-1$ then $\depth M_{i,-1}\geq 2n-i$, so  $\pd M_{i,-1}\leq i$. Thus the assertion of linearity from step $i+1$ on is trivial. On the other hand, for $i\geq g$ we have $\depth M_{i,-1} \geq  2n-3-i$, hence $\pd M_{i,-1} \leq i+3$ and again linearity is trivial. Finally, if $j\geq 0$
then by \ite{linearity} the resolutions of $M_{i-1, j-1}$ and $M_{i-1,j}$ become linear at least one step earlier than what is asserted 
for the resolution of $M_{i,j}$. Now the exact sequence $(*_{i-1,j-1})$ of \ite{hooks} and the long exact sequence for Tor prove \ite{linearity} for $M_{i,j}$.

\medbreak

This completes the proof by induction on $i$.

\medbreak
Having proven that $\depth M_{i,j} = \depth (IR_i)^j$ is at least the value given in  Theorems \ref{new}\ref{CM,Bb} and \ref{main}\ref{depth}, we now prove equality. We will handle five regions with different methods of proof.

\begin{itemize}
 \item[(A)] $j\geq 2$ and $i\leq \ell - 5$ or $j\geq i-g$ and $\ell-4\leq i\leq \ell-2$.
\item[(B)] $j\leq \min\{1,i-g-1\}$ and $i\leq \ell-5$.
\item[(C)] $i-g\leq j \leq 1$.
\item[(D)] $j\leq i-g-1$ and $ i\geq\ell-4$.
\item[(E)] $j \geq i-g$ and $i = \ell-1$.
\end{itemize}

Let $h$ be the Hilbert function of the coordinate ring of the Grassmannian of 2-dimensional subspaces of $k^n$, and write
$\Delta$ for the difference operator, so that $\Delta(h)(t) = h(t) -h(t-1)$. 

To deal with region (A) we need a result that holds in a slightly
larger region:

\begin{lemma}\label{6}\label{socle}
With hypotheses as above, for $(i,j)$ such that
$j\geq -1$ and $i\leq \ell - 5$ or $j\geq i-g$ and $\ell-4\leq i\leq \ell-2$,
the socle of $\rH^4_\mm(M_{i,j})$ is concentrated in degree $-4$ and has $k$-dimension $\Delta^i(h)(j-2)$. 
\end{lemma} 

\begin{proof}[Proof of Lemma~\ref{socle}] 
The concentration statement follows from the fact that the module $M_{i,j}$ has depth $\geq 4$  by \ite{depth for M} and its resolution 
is linear in position $2n-4$ by \ite{linearity}. For the same reason, the dimension of the socle of $\rH^4_\mm(M_{i,j})$ is equal to the 
 $(2n-4)$-th Betti number of $M_{i,j}$. 

We prove by induction on $i$ that this Betti number is $\Delta^i(h)(j-2)$. Suppose that $i=0$. In this case $M_{i,j}=I^j$ and the assertion  
follows from Lemma~\ref{Raicu} if $j\geq 2$ and from the inequality $\depth M_{i,j} \geq 5$ if $j\leq 1$.

Next let $i >0$. If $j=-1$ then $M_{i,j}$ has depth at least  $5$ by \ite{depth for M} and so the $(2n-4)$-th Betti number is $0$, as claimed.
  If $j\geq0$ we use the exact sequence from \ite{hooks},
 $$
 \qquad 0\rTo M_{i-1,j-1} \rTo M_{i-1,j} \rTo M_{i,j} \rTo 0 \, .
 $$
The resolutions of the modules in this sequence have length at most $2n-4$ and are linear in position $2n-4$. Moreover by \ite{linearity},
the resolution of the leftmost module is linear in position $2n-5$. Now it follows from the longexact sequence of Tor that $(2n-4)$-th Betti numbers 
are additive on this sequence, as claimed.\end{proof}

 To deal with region (B) we will use
 
\begin{lemma} \label{generators} For $(i,j)$ in the region (B), the module $\rH^{2n-3-i}_\mm(M_{i,j})$ is generated in degree $0$ and has minimal number of generators
$\Delta^i(h)(i-g-1-j)$.
\end{lemma}
\begin{proof}[Proof of Lemma~\ref{generators}]
We first observe that $i \geq g-1$ and so $j\leq i-g+2$. Hence the duality isomorphism of Theorem~\ref{main}\ref{loco duality}(2)
applies and gives
$$
\rH_\mm^{2n-3-i}(M_{i,j}) \cong \rH_\mm^4(M_{i, i-g+1-j})^\vee(4).
$$
Note that $(i,i-g+1-j)$ is in region (A), giving the desired result through Lemma~\ref{socle}.
\end{proof}

We can now prove the depth equalities in all the regions.
For region (A) we must show that $\rH^4_\mm(M_{i,j}) \neq 0$. Because the coordinate ring of the Grassmannian is Cohen-Macaulay of dimension $\ell=2n-3$, the value of $\Delta^i(h)(t)$ for $i\leq \ell-1$ is the value of the Hilbert function of a ring of dimension $\ell -i >0$ at $t$, which is nonzero
for all $t\geq 0$. In region (A) it is clear that $i\leq \ell-2$ and $j-2\geq 0$, so Lemma~\ite{socle} concludes the proof.

For region (B) we must show that $\rH^{2n-3-i}_\mm(M_{i,j}) \neq 0$. Because $i-g-1-j\geq 0$, we may argue as above using  
Lemma~\ref{generators}.

In region (C) the asserted depths are equal to the dimensions, except when $j=1$ and $i \leq g-2$. In this case $R_i$ has depth
$2n-i \geq n+3$ and $R_i/IR_i=R/I$ has depth $n+1$, showing that $M_{i,1}=(IR_i)(2)$ has depth $n+2$, as claimed.

For region (D) we do decreasing induction on $i$. For $i = \ell-1$ we have, by \ite{hooks} and \ite{iso of M}, an inclusion
 $$
 (IR_{\ell-1})^j(2j) \cong M_{\ell -1,j}\rTo^{a_\ell} M_{\ell-1,j+1}\cong (IR_{\ell-1})^{j+1}(2(j+1)).
 $$
 Because $J = (a_1,\dots, a_\ell)$ is a reduction of $I$, and $I$ is locally a complete intersection on the punctured spectrum, it follows
that $I=J$  on the punctured spectrum, and thus the cokernel of this inclusion has finite length. We know that $M_{\ell-1,j+1}$ has depth $\geq 1$ by \ref{depth for M}. Suppose that
 $M_{\ell-1,j}$ has depth $\geq 2$. It would follow that the inclusion above is an
 isomorphism. If $j=-1$ then 
 $$
 M_{\ell-1,-1}(2) = (IR_{\ell-1})^{-1} = (a_\ell R_{\ell-1})^{-1}\cong R_{\ell-1}(4),
 $$
 so $R_{\ell-1}  = M_{\ell-1, 0}$ would also have depth $\geq 2$, and it suffices to 
 prove the depth equalities for $j\geq 0$.

In this case and still assuming $\depth M_{\ell-1, j} \geq 2$, we would have, by the isomorphism above,
\begin{align*}
I^{j+1} &= a_\ell I^{j}+(J_{\ell-1}:I)\cap I^{j+1}\\
&=a_\ell I^j+J_{\ell-1} I^j \quad \hbox{by $(**)$}\\
&=J I^j \, .
\end{align*}
This implies that the reduction number of $I$ is $\leq (\ell-1)-g -1 = n-4$; however the reduction number of $I$ is $n-3$ (\cite[Proposition 4.2]{EHU}), a contradiction. 
Thus $\depth M_{\ell-1, j} = 1$ in region (D).

For the induction step we use the exact sequence $(*_{i,j})$. Since $\depth M_{i,j+1} \geq 2n-3-i$ by \ite{depth for M} and $\depth M_{i+1,j+1} = 2n-3-i-1$ by the induction hypothesis,
we see that $\depth M_{i,j} = 2n-3-i$ as required.

For $(i,j)$ in region (E) we remark that $\dim M_{i,j} = 4$ so certainly $\depth M_{i,j} \leq 4$, concluding the proof of equality.

\def\rExt{{\rm Ext}}

\medskip 
We now show that for  $n\geq 5$ the ring $R_{g+1}$ is Buchsbaum, which is a statement in Theorem~\ref{new}(a). Although our proof does not apply, we have verified by computation
that it is also Buchsbaum in the case $n=4$ in an example.

By Theorem~\ref{main}\ref{loCo} this ring has only one nonzero intermediate cohomology module, which is 
$\rH^{n-3}_\mm(R_{g+1})$ by the parts of Theorem~\ref{new}(a) that we have already proven. We must show that this module is annihilated by $\mm$.

By the duality statement of Theorem~\ref{main}\ref{loco duality} equation (2) it suffices to prove that
$\rH^{4}_\mm((IR_{g+1})^2) \cong \rH^{n-3}_\mm(R_{g+1})^\vee$ is annihilated by $\mm$.
We will show that this module is $k(1)$.

By Theorem~\ref{main}\ref{depth} we have $\depth IR_i\geq n+1 \geq 6$ for all
$i \leq g$ and thus from the exact sequences $(*_{0,1}),\dots, (*_{g,1})$ of \ite{hooks}
we get $\rH_\mm^4(I^2) \cong \rH_\mm^4((IR_{g+1})^2)$. By \ite{depth for M} and \ite{linearity} 
the resolution of $I^2$ is linear of length $2n-4$, and by Lemma~\ref{Raicu} the last Betti number of this resolution is $1$.  So
$\rExt_S^{2n-4}(I^2,\omega_S)$ is linearly presented and generated by one element of degree one . Since $I$ is a complete intersection on the punctured spectrum, this module has finite length, and thus
must be $k(-1)$ as claimed.
This completes the proof of the Buchsbaum property and hence the proof of Theorem~\ref{new}(a).

\bigbreak
For the regularity statement Theorem~\ref{new}\ref{regularity}, first note that
for $i\leq g-1$ the ring $R_i$ is a complete intersection of $i$ quadrics, and thus
has regularity $i$. If $i \geq g-1,$ then by Theorem~\ref{new}\ref{can} the a-invariant 
satisfies $a(R_i)=-4.$ For $i=g$ the ring $R_i$ is Cohen-Macaulay and so
its regularity is $a(R_i)+ \dim R_g
=i-2.$  If $i=g+1 \geq 5,$ then by Theorem~\ref{new}\ref{unmixed} the ring $R_i$ has only one intermediate local cohomology,
which has been computed above as $\rH_\mm^{n-3}(R_{g+1}) = k(1). $ Thus the regularity is 
$\max\{-1 + n-3, -4+\dim R_{g+1}\} = n-4 = i-4$. This completes the proof of Theorem~\ref{new}\ref{regularity}.

\bigbreak
Now we prove Theorem~\ref{new}\ref{unmixed}. The ideal $IR_i$ contains a nonzerodivisor
on $R_i$ by Proposition~\ref{known facts}\ref{nzd} and is $S_2$ as a module by Theorem~\ref{main}\ref{loco duality}.
From this it follows that $IR_i$ is unmixed of codimension 1.

To prove the two depth inequalities, first note that by Theorem~\ref{main}\ref{depth}
the depth of the modules $IR_i$ and $(IR_i)^{-1}$ is greater or equal to the depth of $R_i$. Thus it suffices
to show that the natural maps
$$
\rH_\mm^{\depth R_i} (IR_i) \to \rH_\mm^{\depth R_i} (R_i) \to \rH_\mm^{\depth R_i} ((IR_i)^{-1})
$$
are both injective. 
Recall that $2 \leq \depth R_i  = 2n-3-i \leq \dim R_i-1$ by Theorem~\ref{new}(a).
Using
the duality isomorphisms in Theorem~\ref{main}(d)(2)
we see that it is enough to prove
the surjectivity of the natural maps in the bottom row of the following commutative diagram:
$$
\begin{diagram}
 \rH_\mm^{4} (I^{i-g})& \lTo& \rH_\mm^{4} (I^{i-g+1}) &\lTo& \rH_\mm^4 (I^{i-g+2})\\
 \dTo&&\dTo&&\dTo\\
 \rH_\mm^{4} ((IR_i)^{i-g})& \lTo& \rH_\mm^{4} ((IR_i)^{i-g+1}) &\lTo& \rH_\mm^4 ((IR_i)^{i-g+2})\, .
\end{diagram}
$$

\medskip
\noindent
By local duality the surjectivity of the maps in the top row follows from the injectivity of the natural maps
$\rExt^{2n-4}(I^j, S) \to \rExt^{2n-4}(I^{j+1}, S)$, which is stated in
\cite[Display (4.8) in the proof of Theorem 4.5]{RWW}  for $j \geq2$ and is obvious for $j\leq 1$.

We claim that the two left hand vertical maps are both surjective, which will complete the argument. This follows using the sequences $(*_{k,j})$ of \ite{hooks} for $0\leq k\leq i-1$
and $i-g-1\leq j\leq i-g$, together with the vanishing of $\rH^5_\mm(M_{k,j})$
in that range of $k,j$ from \ite{local coh}, which applies because
$$
j\geq i-g-1\geq k-g
$$
 and 
$$
\min\{n-j+2, 2n-k-1\}\geq \min\{n-i+g+2, 2n-i\}=2n-i \geq 5 \, .
$$
This completes the proof of Theorem~\ref{new}\ref{unmixed}, and with it all of Theorems~\ref{main} and~\ref{new}, as well as items 
\ite{local coh}--\ite{linearity}.
\end{proof}

We are indebted to Claudiu Raicu for providing a proof of the next result.

\begin{lemma}
\label{Raicu}
If $j\geq 2$ then the last Betti number of $I^j$ is the minimal number of generators of $I^{j-2}$.
\end{lemma}

\begin{proof}
Minimal free resolutions of $I^j$ are computed in~\cite[Theorem 5.4]{ABW} and a description of the modules in the resolutions can be
obtained from ~\cite[Corollary 4.13]{ABW}. For a more explicit version of these resolutions and their modules see also \cite[Theorem 3.1]{RW} with $a=2$ and
\cite[Formulas (1.6), (1.4), (1.5)]{RW}. 
\end{proof}

\begin{lemma}\label{linearitylemma}
Let $N\subset M$ be finitely generated graded modules over a Noetherian positively graded ring $S$ with $S_{0}$ a field. Write $G_{\bullet}$, $F_{\bullet}$ for the minimal graded free resolutions of $N$, $M$ respectively. Suppose that $M$ is generated in degrees $\geq 0$. 
\begin{enumerate}
 \item[$($a$)$] If, for some integer $q\geq \pd M$, the module 
 $G_{q}$ is generated in degree $q$, then $\pd (M/N) \leq q$.
\item[$($b$)$] If  in addition $G_{q-1}$ is generated in degree $q-1 ,$ then
 $$
 \Ext^q(M,S)_{-q} \cong \Ext^q(N,S)_{-q}\oplus \Ext^q(M/N,S)_{-q}.
 $$
  \end{enumerate}
\end{lemma}
\begin{proof}
Consider the diagram of a map of complexes $\phi$ corresponding to the inclusion $N\subset M$:
$$
\begin{diagram}[small]
\cdots &\rTo &G_q &\rTo & G_{q-1}&\rTo&\cdots\\
&&\dTo&&\dTo\\
0&\rTo &F_q &\rTo & F_{q-1}&\rTo&\cdots\\
\end{diagram}
$$

\medskip

\noindent
For item (a), note that the mapping cone of $\phi$ is a possibly nonminimal resolution of $M/N$, which is generated in non-negative degrees.
The module $G_q$ is  the $(q+1)$-st term in this non-minimal resolution, but is generated in degree $q$; however the minimal 
generators of the $(q+1)$-st module in the minimal resolution of $M/N$ must have degrees $\geq q+1$, so $G_q = 0$.
For item (b) we note in addition that the map $G_{q-1}\to F_{q-1}$  must be the inclusion of a direct summand. Thus the $q$-th (last) term in the 
minimal resolution of $M/N$ is $F_q/G_q$.
\end{proof}

\begin{remark}
  If $G_{q-1}$ is generated in degrees $q-1, q$ and $F_{q}$ is generated in degree $q$,
 then the $q$-th module in the minimal graded free resolution of $M/N$ is generated in degree $q$.
\end{remark}

\begin{lemma}\label{duality lemma}
 Let $R$ be a positively graded Noetherian ring over a field $k$, with maximal homogeneous ideal $\m$, and let $M$ be a 
 finitely generated graded $R$-module of dimension $d$ that is Cohen-Macaulay locally on the punctured spectrum and
 equidimensional. For $2 \leq p \leq d-1$ there are isomorphisms of graded $R$-modules
 $$
 \rH_\mm^p(\Hom(M, \omega_R)) \cong \rH^{d+1-p}_\mm(M) ^{\vee} \, ,
$$
where $-^{\vee}$ denotes $k$-dual.
\end{lemma}
\begin{proof} Mapping a graded polynomial ring over $k$ onto $R$ and factoring out an ideal generated by a homogeneous regular sequence 
of maximal length contained in the annihilator of $M$, we may assume that $R$ is Cohen-Macaulay of dimension $d$. Note that the 
local cohomology modules in question have finite length by our assumption on $M$, so the $k$-dual coincides with the graded $k$-dual. Thus by graded local 
duality, it suffices to show that
 $$
 \rH_\mm^p(\Hom(M, \omega_R)) \cong \Ext_R^{p-1}(M,\omega_R).
 $$
  Let $\FF_{\bullet}$ be a graded free resolution of $M$ over $R$.
 Using the fact that $\Hom(\FF_{\bullet}, \omega_R)$ has finite length cohomology except in cohomological degree 0, we see from the 
 exact sequences in local cohomology of the boundaries and cycles that 
 $$
 \rH_{\mm}^p(\Hom(M, \omega_R))  \cong \rH^0_\mm(\Ext_R^{p-1}(M,\omega_R)).
 $$
 Since $\rH^0_\mm(\Ext_R^{p-1}(M,\omega_R)) \cong \Ext_R^{p-1}(M,\omega_R)$,
  the assertion follows.
\end{proof}

\section{The class group of $R_i$}

One reason for interest in the powers $(IR_i)^j$ is that, at least in the generic case,
they constitute the class group of $R_i$. This is the content of the next theorem, a version of 
which was proved in \cite[Theorem 3.4]{HU2} under more stringent depth hypotheses on 
the Koszul homology of $I$.

\begin{theorem}\label{class group}
Let $$
 X = \begin {pmatrix} x_{1,1}&\dots& x_{1,n}\\
  x_{2,1}&\dots&x_{2,n}
  \end {pmatrix}
  $$
  be a generic $2\times n$ matrix over the ring 
 $S = k(\{z_{p,q,r}\})[x_{1,1},\dots,x_{2,n}]$ with $n\geq 4$, where $1\leq p\leq \ell:=2n-3$, and $1\leq q<r\leq n$. 
 Let $\delta_{q,r}$ be the minor of $X$ involving columns $q,r$, and let $I= I_2(X)S$ be the ideal generated by these minors.
 Let $a_p = \sum_{q,r} z_{p,q,r}\delta_{q,r}$, be a generic linear combination of the minors and set
 $J_i = (a_1,\dots, a_i)\subset I,\ K_i =  J_i:I$, and $R_i = S/K_i$, which is a geometric $i$-residual intersection of $I$.

\begin{enumerate}
\item[$($a$)$] If $g-1=n-2\leq i\leq \ell -2$ then $R_i$ is a normal domain, nonsingular in codimension 
 $$
 \min\{\dim R_i-1,  i-n+4\}.
 $$
 \item[$($b$)$]  The divisor class group of $R_i$ is generated by the class of $IR_i$.
\item[$($c$)$]  For positive $j$ the power $(IR_i)^j$ is equal to the symbolic power $(IR_i)^{(j)}$ and $(IR_i)^{-j}=(K_{i+1}R_i)^{(j)}(2j)$.
\end{enumerate}

\end{theorem}
\begin{proof} Recall that $R_i$ is equidimensional of codimension $i$ by Corollary~\ref{residually S2}, and notice that for 
every $q \in V(K_i) \setminus V(IS)$ the localization $(R_i)_q=S_q/(a_1, \ldots, a_i)$ is regular by the generic choice of 
$a_1, \ldots, a_i$. 
The ring $S/I$ has an isolated singularity, so after localizing
 at any prime of $S$ other than the maximal homogeneous ideal, 
 item (a) (except for the easy case $i =g-1$)
 follows from \cite[Lemma 2.3]{HU2} and item (b) as well as (c) for $j>0$ are a consequence of \cite[Theorem 3.4]{HU2}.
   
 It follows from (the statement or the proof of) Theorem~\ref{main}\ref{depth} that the rings $R_i$ and the powers $(IR_i)^j$ have 
 depth $\geq 2,$ and thus the maximal ideal of $R_i$ is not an associated prime of $R_i/(IR_i)^j ,$ proving that $R_i$ is normal
 and
 $(IR_i)^j = (IR_i)^{(j)}$ for $j>0.$

To complete the prove of (c) we notice that $(K_{i+1}R_i)^{(j)} =a_{i+1}^j(IR_i)^{-j} \cong (IR_i)^{-j}(-2j).$
The equality holds because $a_{i+1}$  is a nonzerodivisor  on the domain $R_i,$ the ideals 
on both sides are unmixed and have codimension 1 by Proposition~\ref{known facts}\ref{unmixed}, and $IR_i+K_{i+1}R_i$ has codimension $\geq 2$ as the residual intersection
$K_{i+1}$ is geometric.

To prove item (b) we first show that $K_{i+1}R_i$ and $IR_i$ are prime and $IR_i$ has codimension $1$. As to the ideal $IR_i,$
note that it is prime of codimension $1$ locally on the punctured spectrum by \cite[Theorem 3.4(a)]{HU2}.
So it suffices to prove that the punctured spectrum of $R/IR_i$ is connected and the maximal homogeneous 
ideal is not an associated prime of this ring.
For both these, we only need to show that $R_i/(IR_i)$ as depth $\geq 2$, which follows from
Theorem~\ref{new}\ref{unmixed}.

For the case of $K_{i+1}R_i$, 
we wish to show that $R_{i+1}\cong R_i/K_{i+1}R_i$ is a domain.
Since $a:=a_{i+2}$ is a nonzerodivisor on $R_{i+1}$ by Proposition \ref{known facts}\ref{nzd}, 
it suffices to show that 
$R_{i+1}$ becomes a domain after inverting $a$. This holds because $(R_{i+1})_{a}=S_{a}/ (a_1, \ldots, a_{i+1})$
is a ring of fractions of a polynomial ring over $k[x_{1,1}, \ldots, x_{2,n}]$ (for details see \cite[Lemma 2.2 and the proof of Theorem 3.4(a)]{HU2}).
 
 Also, we know that
 $K_{i+1}\cap IR_i = (a_{i+1}R_i)$ because $R_i$ is a normal domain and $IR_i+K_{i+1}R_i$ has codimension $\geq 2.$ Now item (b) 
 follows from Nagata's Lemma since both $IR_i$ and $K_{i+1}R_i$ are prime ideals, $IR_i$ has codimension $1,$ and the ring $(R_i)_{a_{i+1}}$ is
 factorial as a ring of fractions of a polynomial ring over $k.$
\end{proof}
\vspace{.2cm}

\section{The asymptotic depth of $(IR_i)^j$}\label{asymp sect}

We will show that the maximal Cohen-Macaulay property when $i=\ell-1$ mentioned after Theorem~\ref{main} is not possible for $s$-residual intersections with $s\neq \ell-1$.

\begin{proposition}
 If $A$ is a standard graded Noetherian ring over a local ring $S$, and $M$ is a finitely generated graded $A$-module, then $\depth M_j$ is constant for $j\gg 0 .$
\end{proposition}

\noindent
This proposition is part of~\cite[Theorem 3.3]{CJR}; we give an elementary proof.

\begin{proof}
 We may assume that the residue field of $S$ is infinite, and we proceed by induction on the minimal number of generators $\mu_S(A_1)$. If $\mu_S (A_1)=0$ then $M_j = 0$ for
 $j \gg 0$. 
 
Now suppose that $\mu_S(A_1)>0$, and that the result fails for $M$. Choose $x\in A_1$ general. Truncating $M$, we may assume that $0:_M x = 0$, and we let $N = M/xM$, so that we have an exact sequence
$$
0\rTo M(-1) \rTo^x M\rTo N \rTo 0 \, .
$$
By induction, the proposition holds for $N$, so $\depth N_j$ is constant for $j\geq n\gg 0$.

Let $t$ be minimal so that  $\depth M_j = t$ for infinitely  many $j$. It follows that there exists $j>n$ with $\depth M_j = t < \depth M_{j-1}$. From
the exact sequence above we see that $\depth N_j= t$. Thus the stable value of the depths of the components of $N$ is $t$.

Let $r$ be the second smallest value of $\depth M_j$ that is taken on infinitely many times, and observe that there are then infinitely many values
$j\geq n$ such that $\depth M_j = r \neq \depth M_{j-1}$. Since $\depth N_j = t<r$, it follows that $\depth M_{j-1}$
$=t+1$. From the definition of $r$,
we get $r=t+1$, contradicting the assumption that $\depth M_{j-1} \neq r$.
\end{proof}

\begin{corollary}\label{asymptotic}
 If $R$ is a Noetherian local ring and $I\subset R$ an ideal of positive codimension, then
$$
\ell(I) \leq \dim S - \limsup \depth I^j +1.
$$
\end{corollary}

\begin{proof}
 Applying the proposition to the Rees algebra $\RR$ of $I$ we see that $\depth I^j$ is constant for $j\gg 0$, and replacing
$I$ with a power of $I$ we may assume that $t:=\depth I = \depth I^j$ for all $j\geq 1$. 

Let $\gm\subset S$ be the maximal ideal. It follows that 
$\grade(\gm, \RR) = t$, so $\codim \gm\RR\geq t$
and therefore
$$
\ell(I) = \dim \RR/\gm\RR\leq \dim \RR -\codim \gm \RR \leq (\dim S+1) -t.
$$
\end{proof}

\begin{corollary}
Let $S$ be a Noetherian local ring with infinite residue field, and let $J \subset I\subset S$ be ideals. Suppose that $J$ is generated by $s$ general
elements of $I$ and set $R:= S/(J:I^\infty).$ 
If $R \neq 0$ and $(IR)^j$ is a maximal Cohen-Macaulay $R$-module for infinitely many $j$, then $s = \ell(I)-1$.
\end{corollary}

\begin{proof}
Notice that $s\leq \ell(I)-1$ because otherwise $R =0$. Since $\grade \,  IR >0$ it follows from Corollary 3.2 that $\ell(IR) \leq 1.$
By the Artin-Rees Lemma $(J:I^{\infty}) \cap I^j= J \cap I^j$ for $j \gg 0$ and so the large powers 
of $IR$ and of $I(S/J)$ coincide. This shows that $\ell(IR) = \ell(I(S/J))$. On the other hand, since $J$ is generated by $s$ general elements of
$I$, \cite[Proposition 8.6.1]{SH} gives $\ell(I) \leq s+ \ell(I(S/J))$, which shows that $\ell(I) \leq s +1$ as required.
\end{proof}

\bigskip
\bigskip
\bigskip

\vbox{\noindent Author Addresses:\par
\smallskip
\noindent{David Eisenbud}\par
\noindent{Mathematical Sciences Research Institute,
Berkeley, CA 94720, USA}\par
\noindent{de@msri.org}\par

\medskip
\noindent{Bernd Ulrich}\par
\noindent{Department of Mathematics, Purdue University, West Lafayette, IN 47907, USA}\par
\noindent{bulrich@purdue.edu}\par

}

\end{document}

%% file: preamble.tex
\def\antiddot{\mathinner{\mkern1mu\raise1pt\vbox{\kern7pt\hbox{.}}\mkern2mu
        \raise4pt\hbox{.}\mkern2mu\raise7pt\hbox{.}\mkern1mu}}

\newcommand{\FF}{{\mathbb F}}

\newcommand{\RR}{{\mathbb R}}

\newcommand{\Ext}{{\rm{Ext}}}

\newcommand{\punkt}{\hspace{-.3ex}\raise.15ex\hbox to1ex{\Huge.}}

\def \fix#1 {{\hfill\break \bf (( #1 ))\hfill\break}}

\DeclareMathOperator{\Sym}{Sym}

\DeclareMathOperator{\Hom}{Hom}

\DeclareMathOperator{\depth}{depth}

\DeclareMathOperator{\pd}{pd}

\DeclareMathOperator{\codim}{codim}

\newcommand{\gm}{\mathfrak m}

\newtheorem{theorem}{Theorem}[section]
\newtheorem{lemma}[theorem]{Lemma}
\newtheorem{proposition}[theorem]{Proposition}
\newtheorem{corollary}[theorem]{Corollary}
\newtheorem{conjecture}[theorem]{Conjecture}

\theoremstyle{definition}

\newtheorem{remark}[theorem]{Remark}